\documentclass[final]{article}
\usepackage[latin9]{inputenc}
\usepackage{fullpage}
\usepackage{cite}
\usepackage{amsmath}
\usepackage{url,xspace}
\newcommand{\I}{\ensuremath{\mathbf{I}}\xspace}
\renewcommand{\v}[1]{\boldsymbol{#1}} 
\renewcommand{\i}{\ensuremath{\v{i}}\xspace}
\renewcommand{\j}{\ensuremath{\v{j}}\xspace}
\renewcommand{\k}{\ensuremath{\v{k}}\xspace}
\newcommand{\inner}[2]{\ensuremath{\langle #1, #2 \rangle}\xspace}
\newcommand{\norm}[1]{\ensuremath{||#1||}\xspace}
\let\oldmu\mu\renewcommand{\mu}{\v{\oldmu}} 
\let\oldnu\nu\renewcommand{\nu}{\v{\oldnu}} 
\let\oldxi\xi\renewcommand{\xi}{\v{\oldxi}}
\newcommand{\matlab}{Matlab\textregistered\xspace}
\newcommand{\QTFM}{\textit{Quaternion Toolbox for \matlab}\xspace}
\begin{document}
\title{Fast complexified quaternion Fourier transform}
\author{Salem~Said\footnotemark[2] \and Nicolas~Le~Bihan\footnotemark[2]
   \and Stephen~J.~Sangwine\footnotemark[3] \footnotemark[4]}

\maketitle
\renewcommand{\thefootnote}{\fnsymbol{footnote}}
\footnotetext[2]{Laboratoire des Images et des Signaux, CNRS
                 UMR 5083, ENSIEG, 961 Rue de la Houille Blanche,
                 Domaine Universitaire, BP 46, 38402 Saint Martin d'Hères,
                 Cedex, France. \texttt{nicolas.le-bihan@lis.inpg.fr}}
\footnotetext[3]{Department of Electronic Systems Engineering, University of Essex,
                 Wivenhoe Park, Colchester CO4 3SQ, United Kingdom. \texttt{s.sangwine@ieee.org}}
\footnotetext[4]{The work presented here was carried out at the Laboratoire des Images et des Signaux,
                 Grenoble. Financial support from the Royal Academy of Engineering of the United Kingdom
                 and the Centre National de la Recherche Scientifique (CNRS) in France, 
                 is gratefully acknowledged.}

\begin{abstract}
A discrete complexified quaternion Fourier transform is introduced. This is a generalization of
the discrete quaternion Fourier transform to the case where either or both of the signal/image and
the transform kernel are complex quaternion-valued. It is shown how to compute the transform using
four standard complex Fourier transforms and the properties of the transform are briefly discussed.
\end{abstract}

\section{Introduction}
Discrete quaternion Fourier transforms have been described by several authors
\cite{Sangwine:1996,SangwineEll:1998,BulowSommer:2001,Felsberg:1999,
PeiDingChang: 2001}. The pioneering work of Ell \cite{Ell:thesis} predated all of these
papers and was the inspiration for \cite{Sangwine:1996,SangwineEll:1998}.
Quaternion, or at least vector-valued, signals arise in color image
processing and in processing of signals captured by vector sensors, for example, vector geophones
\cite{BihanMars:2004}. The main motivation for using quaternion algebra to handle such signals is
the correspondence between vector samples and the vector parts of quaternions, which permits
holistic processing of the vector samples and signals, making use of information about the direction
and magnitude of the samples in sample space.

Since the introduction of the quaternion Fourier transform there has been a concern that
perhaps a quaternion-valued signal or image requires a transform based on a `higher' algebra, just
as the Fourier transformation of real-valued signals or images requires a complex transform.
Indeed Bülow and Sommer \cite{BulowSommer:2001} have shown the advantages of using a quaternion
transform even for real (grey-level) images, in order to analyse symmetries in the image. Unlike
Bülow and Sommer, we are interested in signals and images with vector samples (that is, three
components at least, perhaps even three \emph{complex} components).

However, until now there have been no reported definitions of Fourier transforms for
quaternion-valued signals or images based on an algebra of higher dimension than the real
quaternions\footnote{Fourier transforms can be \emph{defined} in many abstract algebras, but
actually computing them numerically by a practical fast algorithm is a different matter.}.

At various times the octonions have been suggested as the basis for a higher transform, but their
algebra is non-associative and therefore presents formidable problems compared to the quaternion
algebra which is of course non-commutative, but otherwise relatively straightforward to work with.
Therefore we have considered the case of a transform based on complexified quaternions, that is
quaternions with four complex components. These can be considered equivalently as complex numbers
with quaternion real and imaginary parts. The properties of this and other quaternion algebras have
been extensively studied in the abstract algebra community, but much less knowledge is available at
the practical algebraic level needed for engineering use. The complexified quaternions were studied
by Hamilton himself \cite{Hamiltonpapers:V3:35} and \cite[Lecture VII, §\,669
(p.664)]{Hamilton:1853}, but he called them \emph{biquaternions}. We prefer to avoid this term as it
was also used by Clifford with a different meaning.

In this paper we introduce the bare basics of complex quaternions necessary to the paper and we
define a complex quaternion Fourier transform analogous to the quaternion Fourier transforms
discussed in \cite{SangwineEll:1998}\footnote{A more recent and more detailed paper by Ell and
Sangwine is under review.}. We then show how the transform may be computed. In common with real
quaternion Fourier transforms, a direct (DFT or FFT) implementation is possible using a complex
quaternion arithmetic package but it is also possible to factorize the transform into four complex
transforms, and thus compute it using standard complex FFT code, thus exploiting the significant
effort expended on performance enhancement of FFT code by other authors.

The details of quaternions and of Fourier transforms in general are outside the scope of this paper.
Ward's book \cite{Ward:1997} provides a good introduction to quaternions, and there are many texts
in the area of signal processing that discuss discrete and fast Fourier transforms. Generalizing the
quaternions to complex quaternions is straightforward -- the four components of the quaternion take
complex values rather than real values, but there are some important ramifications. Firstly,
although the quaternions are a division algebra, the complex quaternions are not. Any real
quaternion (other than zero) has a multiplicative inverse and a non-zero norm or modulus. In
contrast, a complex quaternion (other than zero) may have a zero norm or modulus and therefore no
multiplicative inverse. Obviously this means that care is needed when computing with complex
quaternions, and if we are to define and compute Fourier transforms with them, we must establish
whether vanishing norms will cause problems. A second issue is that a quaternion Fourier transform
necessarily contains as its kernel a quaternion exponential function, and if we are to generalize
this to complex quaternions, we need a complex quaternion exponential function. As has been shown in
\cite{SangwineEll:1998}, quaternion exponentials can be defined using an arbitrary unit pure
quaternion since these are roots of $-1$. (It is fundamental to Fourier transforms that a root of
$-1$ is needed, since the transform analyses a signal into sinusoidal components, and the famous
formula of de Moivre: $e^{i\theta} = \cos\theta + i\sin\theta$ generalizes to any algebra in which a
root of $-1$ can be defined. This follows from the series expansions of the exponential, cosine and
sine functions.) If we are to define a complex quaternion Fourier transform, we require a
\emph{complex} quaternion root of $-1$. We cannot use a real quaternion root of $-1$, or the
transform will simply be a quaternion Fourier transform even if it operates on a complex quaternion
signal (the transform would be equivalent to two quaternion transforms applied separately to the
real and imaginary quaternion parts of the signal). The solutions of the equation
$q^2=-1$~\cite{arXiv:math.RA/0506190} are therefore crucial to defining the transforms presented in
this paper.

\section{Definitions}
Before defining the transform we must define some complexified quaternion concepts. We have not
given proofs of these: most can be found in \cite{Ward:1997} although with different
notation. A complex quaternion has four complex components, for example: $q = w + x\i + y\j + z\k$,
where $w, x, y, z$ are complex and $\i, \j$ and $\k$ are the three quaternion operators. We
denote the complex square root of $-1$ by a capital $\I$ to distinguish it from the quaternion $\i$
with which it must not be confused (this is a fundamental aspect of complexified quaternions).
Note that, although the multiplication of (complex) quaternions is not commutative, $\I$, and
therefore all complex numbers (complex scalars), commute with $\i, \j$ and $\k$, and therefore with
all quaternions. In what follows, $\Re(q)$ denotes the real part of $q$, that is:
$\Re(q) = \Re(w) + \Re(x)\i + \Re(y)\j + \Re(z)\k$, and similarly $\Im(q)$ denotes the imaginary
part of $q$.

\subsection{Inner product, norm and modulus}\label{innerprod}
The inner product of two complex quaternions $q = w_q + x_q\i + y_q\j + z_q\k$ and $p = w_p + x_p\i
+ y_p\j + z_p\k$ is $\inner{q}{p} = w_q w_p + x_q x_p + y_q y_p + z_q z_p$. The semi-norm\footnote{A
semi-norm is a generalization of the concept of a norm, with no requirement that the norm be
zero only at the origin \cite{CollinsDictMaths}. Here, it is possible for $\norm{q}=0$, even though
$q\ne0$ as discussed in subsection~\ref{vanish}.} of a quaternion
$\norm{q} = \inner{q}{q} = w^2 + x^2 + y^2 + z^2$, and the modulus is the square root of the norm.
All of these results are complex, in general. In the special case where $\Im(q)=0$ the semi-norm
reduces to the usual quaternion norm. If the quaternion has zero real part ($\Re(q)=0$), the semi-norm
reduces to a real, but negative, norm (because $\I^2=-1$).
The inner product provides us with the concept of
orthogonality, which is essential to the construction of a (complex) orthonormal basis. Two complex
quaternions are orthogonal if their inner product is zero, that is
$q\perp p \iff \inner{q}{p}=\inner{p}{q}=0$.

\subsection{Complex quaternion roots of $-1$}
Let $\mu$ be a complexified quaternion root of $-1$. It was shown in \cite{arXiv:math.RA/0506190} that
$\mu$ must be a pure complex quaternion satisfying the following conditions:
\[
\Re(\mu)\perp\Im(\mu),\quad \norm{\Re(\mu)} - \norm{\Im(\mu)} = 1
\]
For example, $\mu = \i + \j + \k + (\j - \k)\I$ is a root of $-1$ and has a unit (real) norm,
as can easily be verified.

\subsection{Complex quaternions with null modulus}
\label{vanish}
It is possible for a non-zero complex quaternion to have zero semi-norm. The conditions for the norm
to vanish were discovered by Hamilton \cite[Lecture VII, §\,672 (p.669)]{Hamilton:1853}. Informally,
the result can be shown as follows. Starting with the semi-norm of the quaternion in the form
$\norm{q} = w^2 + x^2 + y^2 + z^2$, expand one of the components in terms of its real and imaginary
parts, for example $w$:
\begin{align*}
w^2 &= (\Re(w) + \Im(w)\I)^2\\
    &= \Re(w)^2 - \Im(w)^2 + 2\Re(w)\Im(w)\I
\end{align*}
We see that the real part of the semi-norm consists of the difference between the norms of the real
and imaginary parts of $q$, and the imaginary part of the semi-norm is twice the inner product of
the real and imaginary parts of $q$. For the semi-norm to vanish, its real and imaginary parts must
separately be zero. This requires the real and imaginary parts of $q$ to be orthogonal real
quaternions of equal norm, that is: $\inner{\Re(q)}{\Im(q)}=0$ and $\norm{\Re(q)}=\norm{\Im(q)}$.

\subsection{Transform pair}
Given the definition of the transform `axis'\footnote{The axis of the transform is a direction in
3-space (in this paper \emph{complex} 3-space) which defines the orientation of the sine component
of the transform in the space of the vector part of the signal samples.} as any complexified
quaternion root of $-1$, $\mu$, we can now define the transform itself, which is almost trivial. In
this paper we restrict our account to the one-dimensional case, but all our definitions and results
generalize to the two-dimensional case without difficulty.
\begin{equation}\label{transform}
\begin{aligned}
F[u] &=            &&\sum_{n=0}^{N-1}\exp\left(-2\pi\mu\frac{nu}{N}\right)f[n]\\
f[n] &= \frac{1}{N}&&\sum_{u=0}^{N-1}\exp\left( 2\pi\mu\frac{nu}{N}\right)F[u]
\end{aligned}
\end{equation}
where the signal $f[n]$ and its `spectrum' $F[u]$ have $N$ samples. The placement of the minus sign
in the forward transform is, of course, arbitrary. The transform pair defined in equation
\ref{transform} must, of course, be evaluated in complex quaternion arithmetic, and the exponential
is a complex quaternion exponential. All the necessary operations, and the transform itself, have
been implemented in the \QTFM library \cite{qtfm}. It is conventional in \matlab for
the scale factor $\frac{1}{N}$ to be applied to the inverse transform and we have followed this
convention in implementing the complex quaternion transforms\footnote{An alternative and widely
used convention is to distribute the scale factor between the forward and inverse transforms. If the
two scale factors are to be equal, each must be $\frac{1}{\sqrt{N}}$}.

An alternative transform is possible, by interchanging the positions of the exponential and the
function $f[n]$. This transform is closely related to the first and easily implemented, although we
omit the details here.

Now we consider whether vanishing norms can occur in the transform. The norm of the exponential
cannot vanish for any $n$ or $u$ because $\mu$ has unit modulus, and the cosine and sine parameters
are real. However, the signal $f[n]$ could contain samples with vanishing norms, unless there is
some restriction on the sample values which prevents the conditions for a vanishing norm occurring.
For example, if the sample values are real or imaginary quaternions they cannot have vanishing norms.
In the more general case where the samples have complex quaternion values, the possibility of
samples with vanishing norm cannot be eliminated. Such samples would not contribute to the `spectrum'
and therefore could not be reconstructed by the inverse transform. It is not known whether samples
of the spectrum could have vanishing norms if computed from real (or imaginary) quaternion signals.
If this did occur, it would again mean that the transform would not correctly invert. Clearly this
point requires further study.

\section{Factorization into four complex transforms}
The definition of the transform in (\ref{transform}) leads directly to a discrete Fourier
transform implementation, which is useful as a reference or test implementation, but useless for
practical computation except for small values of $N$. A fast implementation requires either a
custom coding of a classic FFT algorithm in complexified quaternion arithmetic, or a decomposition
into complex Fourier transforms that can be implemented using existing complex functions or code.
It is the latter approach that we describe here and it follows from a similar approach
developed for the real quaternion case in 2000 \cite{SangwineEll:2000b}.

The signal function $f[n]$ in (\ref{transform}) may be expressed in terms of an orthonormal
basis defined by $\mu$, the transform axis. Apart from being a complex orthonormal basis, this is
exactly as in the real quaternion case. The complex orthonormal basis is defined by $\mu$ and two
other unit pure complex quaternions $\nu$ and $\xi$ such that $\mu\perp\nu\perp\xi$ and $\mu\nu=\xi$
(and $\mu\nu\xi = -1$)\footnote{$\mu\perp\nu$ means $\inner{\mu}{\nu}=0$}. The basis may also be
represented by a $3×3$ \emph{complex} orthogonal (\emph{not Hermitian}) matrix:
\begin{equation*}
\begin{pmatrix}
\oldmu_x & \oldmu_y & \oldmu_z\\
\oldnu_x & \oldnu_y & \oldnu_z\\
\oldxi_x & \oldxi_y & \oldxi_z\\
\end{pmatrix}
\end{equation*}
where $\mu = \oldmu_x\i + \oldmu_y\j + \oldmu_z\k$ \textit{etc}.
We show here only the decomposition of the forward transform, as the inverse is
similar. We first write $f[n]$ in terms of its four complex components:
\[
f[n] = w[n] + x[n]\i + y[n]\j + z[n]\k
\]
Now applying a change of basis, we express $f[n]$ in terms of the new basis. The change of basis is
very simply implemented by resolving the vector part of each sample of $f[n]$ into the three complex
directions defined by the new basis, using the inner product defined in subsection \ref{innerprod}.
Thus we have:
\begin{align*}
x^\prime[n] &= \inner{\mu}{x[n]\i + y[n]\j + z[n]\k}\\
y^\prime[n] &= \inner{\nu}{x[n]\i + y[n]\j + z[n]\k}\\
z^\prime[n] &= \inner{\xi}{x[n]\i + y[n]\j + z[n]\k}
\end{align*}
In the new basis, we have:
\begin{align*}
f[n] &= w[n] + x^\prime[n]\mu + y^\prime[n]\nu + z^\prime[n]\xi\\
     &= \left[w[n] + x^\prime[n]\mu\right] + \left[y^\prime[n] + z^\prime[n]\mu\right]\nu
\end{align*}
Using this change of basis, we are able to separate the transform into the sum of two transforms:
\begin{equation*}
F[u] = \sum_{n=0}^{N-1}\left(
\begin{aligned}
   & e^{-2\pi\mu\frac{nu}{N}}\left[w[n]        + x^\prime[n]\mu\right]\\
+\,& e^{-2\pi\mu\frac{nu}{N}}\left[y^\prime[n] + z^\prime[n]\mu\right]\nu
\end{aligned}
\right)
\end{equation*}
We now separate the terms on the right into real and imaginary parts and group the real
components together and the imaginary components together to make four complex terms.
\begin{equation*}
F[u] = \sum_{n=0}^{N-1}\left(
\begin{aligned}
   & e^{-2\pi\mu\frac{nu}{N}}\left[\Re(w[n]) + \Re(x^\prime[n])\mu\right]\nonumber\\
+\,& e^{-2\pi\mu\frac{nu}{N}}\left[\Im(w[n]) + \Im(x^\prime[n])\mu\right]\nonumber\\
+\,& e^{-2\pi\mu\frac{nu}{N}}\left[\Re(y^\prime[n]) + \Re(z^\prime[n])\mu\right]\nu\nonumber\\
+\,& e^{-2\pi\mu\frac{nu}{N}}\left[\Im(y^\prime[n]) + \Im(z^\prime[n])\mu\right]\nu\nonumber
\end{aligned}
\right)
\end{equation*}
All four of the transforms within this expression are now isomorphic to a complex Fourier transform.
That is, we may replace $\mu$ with $\I$ (the complex root of $-1$, not the quaternion root of $-1$,
$\i$) in order to compute the transform, and we will obtain the same numeric results. After
computing the four complex transforms, all that remains is to re-assemble the parts of $F[u]$ and
invert the change of basis. The latter step is equivalent to multiplying out the factors of $\mu$
and $\nu$ appearing above, but it is more easily performed by a change of basis using the transpose
of the original basis matrix used to change from the standard basis to the $(\mu,\nu,\xi)$ basis.

Note that the factorization into four complex transforms reduces to a factorization into two complex
transforms in the special cases where the function to be transformed has either zero imaginary part
or zero real part. If, additionally, $\mu$ is restricted to be real,
the transform reduces to a previously published \cite{SangwineEll:2000b} quaternion Fourier transform.
However, if the function to be transformed is either real or imaginary, the transform defined in
this paper is not equivalent to a real quaternion transform, because the kernel is a \emph{complex}
quaternion exponential. Finally we note that if the function to be transformed has no vector part
(that is, it is real or complex), and we choose for $\mu$ the degenerate value $\I$ (the complex
root of $-1$) the transform reduces to the usual complex Fourier transform.

The factorization described above has been used in the implementation of the transforms described in
this paper in both one and two dimensions as part of the \QTFM library \cite[v0.8 onwards]{qtfm}. A
direct DFT implementation of equation \ref{transform} is provided for verification of the fast
transforms. The complex fast transforms are the usual \matlab functions, which are implemented
internally by the FFTW library \cite{fftw}.

\section{Discussion}
The complex quaternion transform applied to a real quaternion signal exhibits quaternion conjugate
symmetry in the `spectrum' in the same way as a complex transform applied to a real signal.
The shift and similar properties follow those of the complex transform. It remains to be seen in
what ways the complex quaternion transform may be useful, but one possible area is the development
and fast implementation of linear vector filters based on convolution with complex quaternion
coefficients. We leave a more detailed discussion of these issues for a later paper.

\section{Conclusions}
This paper has introduced for the first time a complex quaternion Fourier transform which may be
used to compute the transform of a quaternion or complex quaternion signal or image. It may be
easily computed using decomposition into four complex Fourier transforms. It appears that the
transform exhibits symmetries that are absent in quaternion Fourier transforms, just as the complex
Fourier transform shows symmetries for real signals that are not apparent with a real-valued
transform such as the Hartley transform or the DCT.

\end{document}